\documentclass[psamsfonts]{amsart}

\usepackage{amssymb,amsfonts}
\usepackage[all,arc]{xy}
\usepackage{enumerate}
\usepackage{mathrsfs}
\usepackage{cite}
\usepackage{enumerate}

\newtheorem{theorem}{Theorem}[section]
\newtheorem{corollary}[theorem]{Corollary}
\newtheorem{proposition}[theorem]{Proposition}
\newtheorem{lemma}[theorem]{Lemma}

\newtheorem{question}[theorem]{Question}
\newtheorem{questions}[theorem]{Questions}

\theoremstyle{definition}
\newtheorem{definition}[theorem]{Definition}

\newtheorem{example}[theorem]{Example}
\newtheorem{examples}[theorem]{Examples}

\newtheorem{remark}[theorem]{Remark}

\DeclareMathOperator{\Frac}{Frac}

\DeclareMathOperator{\Spec}{Spec}
\DeclareMathOperator{\Id}{Id}

\makeatletter

\begin{document}
	
	\title[Dedekind semidomains]{Dedekind semidomains}
	
	\author[P. Nasehpour]{\bfseries Peyman Nasehpour}
	
	\address{Department of Engineering Science \\ Golpayegan University of Technology\\ Golpayegan\\ Iran}
	\email{nasehpour@gut.ac.ir, nasehpour@gmail.com}
	
	
	
	\begin{abstract}
	We define Dedekind semidomains as semirings in which each nonzero fractional ideal is invertible. Then we find some equivalent condition for semirings to being Dedekind. For example, we prove that a Noetherian semidomain is Dedekind if and only if it is multiplication. Then we show that a subtractive Noetherian semidomain is Dedekind if and only if it is a $\pi$-semiring and each of it nonzero prime ideal is invertible. We also show that the maximum number of the generators of each ideal of a subtractive Dedekind semidomain is 2. \end{abstract}
	
	\maketitle
	
	\section{Introduction}\label{sec:intro}
	
	Dulin and Mosher introduce the concept of Dedekind semidomain as a semidomain in which every $k$-ideal is a product of prime $k$-ideals and in addition to some interesting facts, they give the main result of their paper which says that a semisubtractive semidomain $R$ is Dedekind if and only if $R$ is Noetherian, integrally closed, and each of its prime $k$-ideals is maximal \cite[Theorem 1]{DulinMosher1972}. While the statements that they prove in their paper is a nice generalization of some statements for Dedekind domains in the literature, their definition for semidomain is narrow since for them a semidomain is a commutative halfring with multiplicative cancellation and with a multiplicative identity. Note that for them a halfring is a semiring with additive cancellation and a semiring is a halfring if and only if it can be embedded in a ring called its ``ring of differences'' \cite[p. 50]{BrunsGubeladze2009}. Therefore, their paper apparently gives no information on the factorization of ideals into prime ones in those semirings which cannot be embedded in rings.
	
	The main purpose of our paper is to generalize some of the statements for Dedekind domains and prove them for all semidomains (or at least for subtractive semidomains). Before we report what we do in the current paper, we bring some historical information on Dedekind domains which we believe it is useful to explain what our targets in this paper are.   
	
	We may agree that the main result of Dedekind's groundbreaking 1871 work \cite{Dedekind1871} is that every nonzero ideal in the domain of integers of an algebraic number field is a unique product of prime ideals (see Lemma 5.31 in \cite{Jarvis2014} and p. 27 in \cite{Kleiner1998}). Mathematicians have investigated domains having this property and have found many interesting equivalent conditions (definitions) for such domains \cite[p. 143]{Narkiewicz2018}.
	
	Matsumura, in the introduction of his book \cite{Matsumura1989}, says that - a forerunner of the abstract treatment of commutative ring theory - Sono (1886--1969) gives successfully an axiomatic characterization of Dedekind rings in his 1917 paper \cite{Sono1917}. For more on the Japanese mathematician Sono, check \cite{BurtonVanOsdol1995,Sasaki2002}. Then Matsumura adds that Emmy Noether (1882--1935) gives a different system of axioms for Dedekind rings in her 1927 memoir \cite{Noether1927}; the work which is one of the main contributions of Emmy Noether in commutative ring theory as Bourbaki asserts \cite[p. 110]{Bourbaki1994}. 
	
	Dedekind domains (in some resources Dedekind rings \cite[p. 68]{Kaplansky1974} and \cite[p. 82]{Matsumura1989}) have so many equivalent definitions (see Sec. 37 of Chap. VI in \cite{Gilmer1992} and Sec. 4 of Chap. VI in \cite{LarsenMcCarthy1971}). Shortly, we list those of the equivalent statements that either we generalize in this paper or we believe they are useful for the purposes of our paper. For an arbitrary domain $D$, the following statements are equivalent:
	
	\begin{enumerate}
		\item Nonzero fractional ideals of $D$ form a group under multiplication (Krull \cite{Krull1935}).
		
		\item Every nonzero proper ideal of $D$ is a product of prime ideals in $D$ (Matusita \cite{Matusita1944}).
		
		\item The domain $D$ is Noetherian, and for its maximal ideals $\mathfrak{m}$, there are no ideals $\mathfrak{a}$ strictly between $\mathfrak{m}^2$ and $\mathfrak{m}$, i.e., $\mathfrak{m}^2 \subset \mathfrak{a} \subset \mathfrak{m}$ (I.S. Cohen \cite{Cohen1950}).
		
		\item The domain $D$ is Noetherian, and the lattice of its ideals is distributive, i.e., for all ideals $\mathfrak{a}, \mathfrak{b},$ and $\mathfrak{c}$ in $D$, $\mathfrak{a} \cap (\mathfrak{b} + \mathfrak{c}) = (\mathfrak{a} \cap \mathfrak{b}) + (\mathfrak{a} \cap \mathfrak{c})$ (I.S. Cohen \cite{Cohen1950}).
		
		\item The domain $D$ is Noetherian, and $(\mathfrak{a}+\mathfrak{b}) \cdot (\mathfrak{a} \cap \mathfrak{b}) = \mathfrak{a}\mathfrak{b}$, for any two ideals $\mathfrak{a}$ and $\mathfrak{b}$ of $D$ (Jensen \cite{Jensen1963}).
		
		\item Every nonzero ideal in $D$ is an intersection of finitely many powers of prime ideals (Butts and
		Gilmer \cite{ButtsGilmer1966}).
		
		\item The domain $D$ is Noetherian, and if $\mathfrak{a} \subseteq \mathfrak{b}$ are ideals of $D$ with $\mathfrak{a} \mathfrak{b}^n = \mathfrak{b}^{n+1}$ for some $n \geq 1$, then $\mathfrak{a} = \mathfrak{b}$ (Hays \cite{Hays1973}).
	\end{enumerate}
	
	Dedekind domains have played a crucial role in the development of algebraic geometry as well as commutative ring theory (see Dieudonn\'{e} \cite{Dieudonne1972} and Kleiner \cite{Kleiner1995}). For this reason, it is not a surprise to see that some algebraists have attempted to define and investigate the concept of Dedekind domain in other branches of algebra. Dulin and Mosher \cite{DulinMosher1972} investigate Dedekind semidomains. On the other hand, according to the equivalent conditions given by Dorofeeva \cite{Dorofeeva1972}, a Dedekind monoid is a commutative cancellative monoid in which every ideal is a product of prime ideals (also see \cite{DorofeevaMannepalliSatyanarayana1974} and \cite{HannepalliSatyanarayana1974} on Dedekind monoids). All these have motivated the author to investigate this concept in semiring theory, though with a different approach in comparison to the one by Dulin and Mosher.
	
	Since the glossary of the language of semiring theory is not standardized and sometimes confusing \cite{Glazek2002}, before explaining the content of our paper, we need to fix some definitions and terminologies.
	
	In this paper, a semiring is an algebraic structure $(S,+,\cdot,0,1)$ with the following properties:
	
	\begin{enumerate}
		\item $(S,+,0)$ is a commutative monoid,
		\item $(S,\cdot,1)$ is a commutative monoid with $1\neq 0$,
		\item $a(b+c) = ab+ac$ for all $a,b,c\in S$,
		\item $a\cdot 0 = 0\cdot a = 0$ for all $a\in S$. 
	\end{enumerate}
	
	A semiring that is not a ring is called a proper semiring. A semiring $S$ is semidomain if $ab = ac$ implies $b=c$ for all $b,c \in S$ and all nonzero $a\in S$.
	
	Let us recall that a nonempty subset $I$ of a semiring $S$ is said to be an ideal of $S$, if $a+b \in I$ for all $a,b \in I$ and $sa \in I$ for all $s \in S$ and $a \in I$ \cite{Bourne1951}. An ideal $I$ of a semiring $S$ is called proper, if $I \neq S$. A proper ideal $P$ of a semiring $S$ is called a prime ideal of $S$, if $ab\in P$ implies either $a\in P$ or $b\in P$. We collect all prime ideals of $S$ in $\Spec(S)$. An ideal $I$ of a semiring $S$ is called subtractive \cite{Golan1999(b)} (formerly $k$-ideal \cite{Henriksen1958}) if $a+b\in I$ and $a\in I$ imply $b\in I$ for all $a,b\in S$. We say a semiring is subtractive if each of its ideals is subtractive.
	
	An ideal $I$ of a semiring $S$ is finitely generated if it is generated by finitely many elements of $S$. A semiring is Noetherian if each ideal of $S$ is finitely generated. An ideal $I$ of $S$ is principal if it is generated by a single element of $S$. A semidomain is a principal ideal semidomain (for short, PISD) if each of its ideals is principal. For more on principal ideal semidomains, see \cite{Nasehpour2018S}.
	
	A semiring $S$ is a discrete valuation semiring (for short, DVS), if $S$ is a valuation semidomain and its value group is $\mathbb Z$ \cite[Definition 3.1]{Nasehpour2018(b)}. A semiring $S$ is a DVS if and only if $S$ is a PISD possessing a unique maximal ideal \cite[Theorem 3.6]{Nasehpour2018(b)}.
	
	A nonempty subset $W$ of a semiring $S$ is said to be a multiplicatively closed set if $1\in W$ and for all $w_1,w_2 \in W$, we have $w_1 w_2 \in W$. Note that if $W$ is such a set in $S$, one can define the localization of $S$ at $W$, similar to the definition of the localization in ring theory (refer to \cite{Kim1985} and \cite[\S 11]{Golan1999(b)}). Now, it is clear that similar to the concept of the field of fractions in ring theory, one can define the semifield of fractions $F(S)$ of the semidomain $S$ as the localization of $S$ at $S-\{0\}$ \cite[p. 22]{Golan1999(a)}. Similar to ring theory, the localization of $S$ at the multiplicatively closed subset $W=S-\mathfrak{p}$ is denoted by $S_\mathfrak{p}$.
	
	Here is a brief description of the content of the paper: In Section \ref{sec:fracinvert}, we bring the definition of fractional and invertible ideals \cite{GhalandarzadehNasehpourRazavi2017} of a semidomain and prove some facts about them. Note that if $S$ is a semidomain and $K=F(S)$ is its semifield of fractions, a nonempty subset $A$ of $K$ is called a fractional ideal of $S$ if $A$ is an $S$-subsemimodule of $K$ and there exists a nonzero $d\in S$ such that $dA \subseteq S$ (for the concept of semimodule and subsemimodule, see \cite[\S14]{Golan1999(b)}). A fractional ideal $A$ of a semidomain $S$ is invertible if there exists a fractional ideal $B$ of $S$ such that $AB = S$ (see Definition \ref{fractional} and Definition \ref{invertible}). This section mainly includes some preparatory materials for the rest of the paper. For example, in Proposition \ref{FactorizationofInvertible}, we show that if $S$ be a semidomain, then a factorization of an invertible integral ideal $\mathfrak{a}$ of $S$ into prime ideals $\mathfrak{p}_1 \mathfrak{p}_2 \cdots \mathfrak{p}_k$ is unique, up to the permutation of the $\mathfrak{p}_i$s. Note that a fractional ideal of a semiring $S$ is integral if it is a subset of $S$.
	
	Section \ref{sec:DDKsemi} is devoted to Dedekind semidomains. In Definition \ref{Dedekindsemidomaindef}, we define a semidomain to be Dedekind if each of its nonzero fractional ideals is invertible. Note that our definition for Dedekind semidomains is not equivalent to the one given by Dulin and Mosher in \cite{DulinMosher1972} (see Example \ref{DulinMosherEx1} and Example \ref{DulinMosherEx2}). In Theorem \ref{Dedekind1}, we prove that $S$ is a Dedekind semidomain if and only if one of the following statements hold:
	
	\begin{enumerate}
		\item $S$ is a Noetherian Pr\"{u}fer semidomain.
		
		\item $\Frac(S)$ is an abelian group.
		
		\item $S$ is Noetherian and every 2-generated ideal of $S$ is invertible.
		
		\item Each nonzero integral ideal of $S$ is invertible.
	\end{enumerate}

    We add that a semiring $S$ is Pr\"{u}fer if each nonzero finitely generated ideal of $S$ is invertible. By $\Frac(S)$, we mean the set of all nonzero fractional ideals of a semidomain $S$.
	
	Then we continue to find other equivalent definitions for Dedekind semidomains and prove that $S$ is a Dedekind semidomain if and only if one of the following statements hold:
	
	\begin{enumerate}[(a)]
		
		\item $[\mathfrak{a}:\mathfrak{b}]_S +[\mathfrak{b}:\mathfrak{a}]_S=S$ for all ideals $\mathfrak{a}$ and $\mathfrak{b}$ of $S$ (Theorem \ref{Dedekind2}).
		
		\item $S$ is Noetherian and $S_\mathfrak{m}$ is a DVS for each maximal ideal $\mathfrak{m}$ (Theorem \ref{LocalizationofDedekind}).
		
		\item $S$ is a Noetherian multiplication semidomain (Theorem \ref{Dedekindismultiplication}).
		
	\end{enumerate}
	
	Note that in general, we denote $[A:B]_D = \{d\in D : dB \subseteq A\}$, whenever $A$, $B$, and $D$ are sets that for which this makes sense. A semiring $S$ is multiplication  if the condition $\mathfrak{a} \subseteq \mathfrak{b}$ for ideals $\mathfrak{a}$ and $\mathfrak{b}$ of $S$ implies the existence of an ideal $\mathfrak{c}$ satisfying $\mathfrak{a} = \mathfrak{bc}$ (check Definition \ref{multiplicationsemiring}).
	
	Also, similar to ring theory \cite[p. 210]{LarsenMcCarthy1971}, we define a semiring $S$ to be a weak multiplication semiring if $\mathfrak{a} \subseteq \mathfrak{p}$, where $\mathfrak{a}$ is an ideal of $S$ and $\mathfrak{p}$ is a prime ideal of $S$, implies that there exists an ideal $\mathfrak{c}$ of $S$ such that $\mathfrak{a} = \mathfrak{p}\mathfrak{c}$ (see Definition \ref{weakmultiplicationsemiring}). Then in Theorem \ref{weakmultiplicationstrictlybetween}, we prove that if $S$ is a weak multiplication semiring, then for its maximal ideals $\mathfrak{m}$, there are no ideals $\mathfrak{a}$ strictly between $\mathfrak{m}^2$ and $\mathfrak{m}$. This implies that if $S$ is a Dedekind semidomain, then there are no ideals $\mathfrak{a}$ strictly between $\mathfrak{m}^2$ and $\mathfrak{m}$ (see Corollary \ref{Dedekindstrictlybetween}). 
	
	In the final phase of this section, we also obtain some properties for Dedekind semidomains similar to their counterparts in ring theory. For example, in Theorem \ref{Dedekindprimeismaximal}, we prove that each nonzero prime ideal of a Dedekind semidomain is maximal. This statement means that the Krull dimension of a Dedekind semidomain $S$ is at most 1, where similar to ring theory, the Krull dimension $\dim S$ of a semiring $S$ is the supremum of the lengths of all chains of prime ideals in $S$ \cite[p. 578]{AlarconAnderson1994}.
	
	Inspired by the proof of Theorem 12.6.8 in \cite{Steinberger1993}, we show that if $S$ is a semidomain, then $S$ is Dedekind if and only if every nonzero integral ideal of $S$ is uniquely represented as a product of maximal ideals in $S$ (check Theorem \ref{ConverseUFTforIdeals}).
	
	In Theorem \ref{SubtractiveDedekind1}, we find a new characterization for subtractive Dedekind semidomains and we prove that if $S$ is a subtractive semidomain, then $S$ is a Dedekind semidomain if and only if every nonzero prime ideal in $S$ is invertible and $S$ is a $\pi$-semiring. We say $S$ is a $\pi$-semiring if each nonzero proper principal ideal of $S$ is a product of prime ideals of $S$ (Definition \ref{pisemiring}). We also prove that if $S$ is a subtractive Dedekind semidomain, then the number of the generators of each ideal of $S$ is at most 2 (see Theorem \ref{SubtractiveDedekind2}).
	
	Finally in Theorem \ref{FiniteSpecPISD}, we show that a subtractive Dedekind semidomain with only finitely many prime ideals is a PISD.
	
	In Section \ref{sec:Mcancel}, we introduce $M$-cancellation ideals in semirings that is a generalization of the concept of cancellation ideals in semirings introduced by LaGrassa in her Ph.D. dissertation \cite{LaGrassa1995}. Let $S$ be a semiring and $M$ an $S$-semimodule. We define an ideal $\mathfrak{a}$ of $S$ to be $M$-cancellation if for all $S$-subsemimodules $P$ and $Q$ of $M$, $\mathfrak{a}P = \mathfrak{a}Q$ implies $P=Q$ (see Definition \ref{Mcancellationdef}). Invertible ideals of a semiring are some good examples for $M$-cancellation ideals. In addition to other results, in Theorem \ref{MGaussian}, we show that if $S$ is a Pr\"{u}fer (in particular, a Dedekind) semidomain and $M$ is a subtractive $S$-semimodule, then $S$ is $M$-Gaussian. Note that we define a semiring $S$ to be $M$-Gaussian if $c(fg) = c(f)c(g)$ for all $f\in S[X]$ and $g\in M[X]$ (see Definition \ref{MGaussiandef}), where for any $g \in M[X]$, the content of $g$, denoted by $c(g)$, is defined to be the $S$-subsemimodule of $M$ generated by the coefficients of $g$. An $S$-semimodule $M$ is subtractive if each its $S$-subsemimodules is subtractive. As $S$-subsemimodule $N$ of $M$ is subtractive if $x+y\in N$ and $x\in N$ imply $y\in N$, for all $x,y\in M$. At the end, the reader is warned that we use ``$\subseteq$'' for inclusion and ``$\subset$'' for strict inclusion \cite[p. 17]{Monk1969}.
	
	\section{Fractional and Invertible Ideals of Semidomains}\label{sec:fracinvert}
	
	First, we recall the concept of fractional ideals \cite[Definition 1.1]{GhalandarzadehNasehpourRazavi2017}:
	
	\begin{definition}
		
		\label{fractional}
		
		Let $S$ be a semidomain and $K=F(S)$ its semifield of fractions. A nonempty subset $A$ of $K$ is called a fractional ideal of $S$ if the following conditions are satisfied:
		
		\begin{enumerate}
			
			\item $A$ is an $S$-subsemimodule of $K$, that is, if $a, a^{\prime} \in A$ and $ s \in S$, then $a + a^{\prime} \in A $ and $sa \in A$.
			
			\item There exists a nonzero $d\in S$, known as the common denominator of $A$, such that $dA \subseteq S$.
			
		\end{enumerate}
	\end{definition}

	\begin{proposition}
		
		Let $S$ be a semidomain and $K$ its semifield of fractions. Then the following statements hold:
		
		\begin{enumerate}
			
			\item Each ideal of $S$ is fractional.
			
			\item If $K\neq S$, then $K$ is not a fractional ideal of $S$.
			
			\item If $S$ is Noetherian and $I \subseteq K$ is an $S$-semimodule, then $I$ is finitely generated if and only if $dI \subseteq S$ for some nonzero $d\in A$. 
		\end{enumerate}
		
		\begin{proof} (1): The common denominator is $d=1$.
			
			(2): Otherwise, there is a nonzero $d$ in $S$ such that $F = S[d^{-1}]$. This implies that $d^{-2} = d^{-1}s$ for some nonzero $s\in S$ and so, $d^{-1} = s\in S$ showing that $F =  S[d^{-1}] = S$, a contradiction.
			
			(3): For the forward implication, if $x_1 / s_1, x_2 / s_2,\dots, x_n / s_n$ generate $I$ as an $S$-semimodule, then $dI \subseteq S$ for $d=x_1 x_2 \cdots x_n$. Conversely, if $dI \subseteq S$, then $dI$ is an ideal of $S$ and so, it is finitely generated (since $S$ is Noetherian). Suppose $s_1, s_2, \dots, s_n$ generate $dI$. Therefore, $I$ is generated by $s_1 / d, s_2 / d, \dots, s_n / d$.
		\end{proof}
		
	\end{proposition}
	
	\begin{remark}

		Let $S$ be a semidomain and $K$ its semifield of fractions.
		
		\begin{enumerate}
			\item Similar to ring theory, we usually call an ideal of $S$ an integral ideal of $S$. 
			
			\item If $x\in K$, then the cyclic $S$-subsemimodule $Sx$ of $K$ is a fractional ideal of $S$, known as a principal fractional ideal of $S$. We may denote $Sx$ by $(x)$ as well.
			
			\item If $x_1, x_2, \dots, x_n \in K$, then $Sx_1 + Sx_2 + \dots + Sx_n$ is a fractional ideal of $S$, generated by the set $\{x_1,x_2,\dots,x_n\}$. We may denote $Sx_1 + Sx_2 + \dots + Sx_n$ by $(x_1,x_2,\dots,x_n)$ as well. 
			
		\end{enumerate}
	\end{remark}
	
	\begin{examples}
		\begin{enumerate}
			\item Let $\mathbb N_0$ be the set of all non-negative integers. Clearly, the set of all non-negative quotient numbers $\mathbb Q^{\geq 0}$ is its semifield of fractions. Let $n$ be a positive integer. The set $\displaystyle I = \frac{1}{n} \mathbb N_0 = \{\frac{m}{n} : m\in \mathbb N_0\}$ is a fractional ideal of $\mathbb N_0$. Note that $I$ as an $\mathbb N_0$-subsemimodule of $\mathbb{Q}^{\geq 0}$ is generated by $\displaystyle \frac{1}{n} \in \mathbb{Q}^{\geq 0}$ and $nI \subseteq \mathbb N_0$.
			
			\item  Let $I$ be the $\mathbb N_0$-subsemimodule of $\mathbb Q^{\geq 0}$, generated by the quotients of the form $1/2^n$, where $n$ runs over all positive integers. Since there is no positive integer $d$ such that $dI \subseteq \mathbb N_0$, $I$ is not a fractional ideal of $\mathbb N_0$.
		\end{enumerate}
	\end{examples}
	
	\begin{proposition}
		
		\label{betweenprincipal}
		
		Let $S$ be a semidomain and $A$ a nonzero fractional ideal of $S$. Then there are two nonzero elements $c$ and $d$ in $S$ such that $(c) \subseteq A \subseteq (d^{-1})$.
		
		\begin{proof} Since $A$ is nonzero, $A$ contains a nonzero element $y\in K$. It is obvious that some multiple of $y$ is $c\in A \cap S$. Obviously, $(c) \subseteq A$. On the other hand, since $A$ is fractional, there is a nonzero $d\in S$ such that $dA \subseteq S$. This implies that $A \subseteq (d^{-1})$. 
		\end{proof}
	\end{proposition}
	
	Let $S$ be a semidomain, $K=F(S)$ its semifield of fractions, and $A$ and $B$ $S$-subsemimodules of $K$. The sum of $A$ and $B$ is defined as \[A+B := \{a+b : a\in A, b\in B\}.\] Their product is \[A\cdot B := \{\sum^n_{i=1} a_i b_i: a_i \in A, b_i \in B, n\in \mathbb N\}.\] And the residual quotient of $A$ by $B$ is defined as \[[A:B] :=\{x\in K: Bx \subseteq A\}.\]

	\begin{proposition}
		
		\label{sumproductfrac}
		
		Let $S$ be a semidomain and $K=F(S)$ its semifield of fractions. Then the following statements hold:
		
		\begin{enumerate}
			\item If $\{A_1,A_2,\dots,A_n\}$ is a finite set of fractional ideals of $S$, then their sum $\sum_{i=1}^{n} A_i$, their product $\prod_{i=1}^{n} A_i$, and their intersection $\bigcap^n_{i=1} A_i $ are also fractional ideals of $S$.
			\item If $A$ is a fractional ideal of $S$ and $J$ a nonzero $S$-subsemimodule of $K$, then $[A : J]$ is a fractional ideal of $S$.
		\end{enumerate}
		
		\begin{proof} 
			
			The proof of the statement (1) is straightforward. For the proof of (2), it is easy to check that $[A: J]$ is an $S$-semimodule. Let $d\neq 0$ be the common denominator of $A$. Choose a nonzero $t$ in $J \cap S$. Clearly, for any $x\in [A : J]$, $xt \in A$. So, $x(dt)\in S$. Therefore, $dt$ is the common denominator of $[A: J]$ and this completes the proof.
		\end{proof}
		
	\end{proposition}
	
	Now, we recall the concept of invertible fractional ideals \cite[Definition 1.2]{GhalandarzadehNasehpourRazavi2017}:
	
	\begin{definition}
		
		\label{invertible}
		
		A fractional ideal $A$ of a semidomain $S$ is invertible if there exists a fractional ideal $B$ of $S$ such that $AB = S$. If a fractional ideal $A$ is invertible, we denote its inverse by $A^{-1}$.
		
	\end{definition}

	\begin{proposition} 
		
		\label{Dedekindidentity}
		
		Let $S$ be a semidomain. The following statements hold:
		
		\begin{enumerate}
			\item (Dedekind's Identity \cite[p. 113]{Dedekind1895}) For all fractional ideals $A$, $B$, and $C$, the following formula holds: \[(A+B+C)(BC+CA+AB)=(B+C)(C+A)(A+B).\]
			
			\item If each 2-generated fractional ideal of $S$ is invertible, then each nonzero finitely generated fractional ideal of $S$ is invertible.
		\end{enumerate}
		
		\begin{proof} (1): Both sides are $ABC+A^2 B+B^2 A+A^2 C+C^2 A + B^2 C + C^2 B .$
			
			(2): The proof is by induction. Let $n>2$ be a natural number and suppose that all nonzero ideals of $S$ generated by less than $n$ generators are invertible ideals and $L=(a_1,a_2,\ldots,a_{n-1},a_n)$ be an ideal of $S$. If we put $A=(a_1)$, $B=(a_2,\ldots,a_{n-1})$ and $C=(a_n)$, then by induction's hypothesis the ideals $A+B$, $B+C$, and $C+A$ are all invertible ideals. Since the product of fractional ideals of a semiring is invertible if and only if every factor of this product is invertible, by Dedekind's Identity, the ideal $A+B+C=L$ is invertible and the proof is complete. \end{proof}
		
	\end{proposition}
	
	\begin{definition}
		The set of all invertible fractional ideals of $S$ is an abelian group which we call it the ideal group of $S$ and denote it by $\mathcal{I}_{\mathcal{F}}(S)$.
	\end{definition}
	
	\begin{example}
		Let $S$ be a discrete valuation semiring with $(t)$ its maximal ideal. Clearly, its nonzero fractional ideals of $S$ are all in the form $(t^n)$ for $n\in \mathbb Z$, all of which are invertible. Since $(t^m)(t^n) = (t^{m+n})$ for all $m,n\in \mathbb Z$, thus the ideal group of $S$ is isomorphic to $\mathbb Z$ (under addition). For more on valuation semirings, refer to \cite{Nasehpour2018(b)}.
	\end{example}
	
	\begin{proposition}
		
		\label{FactorizationofInvertible}
		
		Let $S$ be a semidomain. A factorization of an invertible integral ideal $\mathfrak{a}$ of $S$ into prime ideals $\mathfrak{p}_1 \mathfrak{p}_2 \cdots \mathfrak{p}_k$ is unique, up to the permutation of the $\mathfrak{p}_i$s.
		
		\begin{proof}
			
			Let $\mathfrak{p}_1 \mathfrak{p}_2 \cdots \mathfrak{p}_k = \mathfrak{q}_1 \mathfrak{q}_2 \cdots \mathfrak{q}_l$ be two factorizations of $\mathfrak{a}$ into prime ideals. Since $\mathfrak{a}$ is invertible, all the $\mathfrak{p}_i$s and $\mathfrak{q}_j$s are invertible. Suppose $\mathfrak{q}_1$ is minimal among  $\mathfrak{q}_1, \mathfrak{q}_2, \dots, \mathfrak{q}_l$. Since $\mathfrak{p}_1 \mathfrak{p}_2 \cdots \mathfrak{p}_k \subseteq \mathfrak{q}_1$, we see that some $\mathfrak{p}_i \subseteq \mathfrak{q}_1$. With the same argument, we can conclude that some $\mathfrak{q}_j \subseteq \mathfrak{p}_i$. By assumption, $\mathfrak{q}_1$ is minimal among the $\mathfrak{q}_j$s. So, $\mathfrak{q}_1 = \mathfrak{p}_i = \mathfrak{q}_j$ and by canceling $\mathfrak{q}_1$, we obtain that $\mathfrak{p}_1 \mathfrak{p}_2 \cdots \mathfrak{p}_{i-1} \mathfrak{p}_{i+1} \cdots \mathfrak{p}_k = \mathfrak{q}_2 \cdots \mathfrak{q}_l$. Therefore, by induction on $k$, we obtain the result.
		\end{proof}
		
	\end{proposition}
	
	\begin{proposition}
		Let $S$ be a semidomain and $A$ and $B$ some fractional ideals of $S$. Then the following statements hold:
		
		\begin{enumerate}
			
			\label{invertibleisfg}
			
			\item $[AB : A]A = AB$.
			
			\item $[S: A]$ is a fractional ideal of $S$.
			
			\item If $A$ is invertible, then $A^{-1} = [S:A]$.
			
			\item If $A$ is an invertible integral ideal of $S$, then $A$ is finitely generated.
		\end{enumerate}
		
		\begin{proof} (1): Straightforward.\\
			(2): Since $S$ is fractional and $A$ an $S$-semimodule, by Proposition \ref{sumproductfrac}, $[S : A]$ is fractional.\\
			(3): In the formula, $[AB : A]A = AB$, put $AB=S$.\\
			(4): Let $A$ be an invertible integral ideal of $S$. So, there is a fractional ideal $B$ of $S$ such that $AB = S$. This implies that $\sum^n_{i=1} x_iy_i=1,$ for some $x_1,x_2,\dots,x_n \in A$ and $y_1,y_2,\dots,y_n \in B$. Clearly, the set $\{x_i\}^n_{i=1}$ generates $A$ in $S$.
		\end{proof}
	\end{proposition}

\begin{proposition}
	\label{reyesfactorizationtheorem}
	
	Let $\mathfrak{a} \subseteq \mathfrak{b}$ be integral ideals of a semiring $S$. If $\mathfrak{b}$ is invertible, then $\mathfrak{a} = \mathfrak{b} \cdot [\mathfrak{a} : \mathfrak{b}]$, where $[\mathfrak{a} : \mathfrak{b}]= \{s\in S: s\mathfrak{b} \subseteq \mathfrak{a}\}.$
	
	\begin{proof}
		
	From $\mathfrak{a} \subseteq \mathfrak{b}$, we infer that $\mathfrak{b}^{-1}\mathfrak{a} \subseteq S$ is an integral ideal of $S$. From $\mathfrak{b}^{-1}\mathfrak{a} \cdot \mathfrak{b} \subseteq \mathfrak{a}$, we obtain that $\mathfrak{b}^{-1}\mathfrak{a} \subseteq [\mathfrak{a}:\mathfrak{b}]$. Since $\mathfrak{b} \cdot \mathfrak{b}^{-1} = S$, multiplying this by $\mathfrak{b}$ we obtain that $\mathfrak{a} \subseteq \mathfrak{b} \cdot [\mathfrak{a} : \mathfrak{b}]$. On the other hand, be definition, $\mathfrak{b} \cdot [\mathfrak{a} : \mathfrak{b}] \subseteq \mathfrak{a}$. This completes the proof. \end{proof}
	
\end{proposition}

	Let us recall that a set of distinct ideals $\{\mathfrak{a}_1, \dots, \mathfrak{a}_n\}$ of $S$ is coprime if $\mathfrak{a}_i + \mathfrak{a}_j = S$ for all $1 \leq i < j \leq n$. Note that the set $\{\mathfrak{m}_1, \dots, \mathfrak{m}_n\}$ of distinct maximal ideals of an arbitrary semiring $S$ is coprime. The proof of the following statement is straightforward and we bring it here only for the sake of reference:
	
	\begin{proposition}
		
		\label{IdealOperationsCoprime}
		
		Let $\{\mathfrak{a}_1, \dots, \mathfrak{a}_n\}$ be a coprime set of ideals of a semiring $S$. Then the following identities hold:
		
		\begin{enumerate}
			\item $(\mathfrak{a}_1^{k_1} \cdots \mathfrak{a}_n^{k_n}) + (\mathfrak{a}_1^{l_1} \cdots \mathfrak{a}_n^{l_n}) = \mathfrak{a}_1^{\min\{k_1,l_1\}} \cdots \mathfrak{a}_n^{\min\{k_n,l_n\}}$.
			
			\item $(\mathfrak{a}_1^{k_1}  \cdots  \mathfrak{a}_n^{k_n}) \cap (\mathfrak{a}_1^{l_1} \cdots \mathfrak{a}_n^{l_n}) = \mathfrak{a}_1^{\max\{k_1,l_1\}} \cdots \mathfrak{a}_n^{\max\{k_n,l_n\}}$
			
			\item $(\mathfrak{a}_1^{k_1}  \cdots  \mathfrak{a}_n^{k_n}) \cdot (\mathfrak{a}_1^{l_1} \cdots \mathfrak{a}_n^{l_n}) = \mathfrak{a}_1^{k_1+l_1} \cdots \mathfrak{a}_n^{k_n+l_n}$
		\end{enumerate}
		
	\end{proposition}

	\section{Dedekind Semidomains}\label{sec:DDKsemi}
	
	According to the equivalent conditions explained on p. 143 in Narkiewicz' book \cite{Narkiewicz2018}, a Dedekind domain is a domain in which nonzero fractional ideals form a group under multiplication; apparently, proved for the first time by Krull on p. 13 in \cite{Krull1935}. Inspired by this, we give the following definition:
	
	\begin{definition}
		
		\label{Dedekindsemidomaindef}
		
		We define a semidomain $S$ to be a Dedekind semidomain if each nonzero fractional ideal of $S$ is invertible.
	\end{definition}
	
	\begin{remark}
		
		For an interesting approach to Dedekind domains see \cite{BerrickKeating2000}. For the roots and applications of Dedekind domains in algebraic number theory refer to \cite{Ash2010,Jarvis2014}. For modern applications of fractional ideals and Dedekind domains see \cite{BrunsGubeladze2009}.	
	\end{remark}
	
	Let us recall that a semidomain $S$ is a Pr\"{u}fer semidomain if each nonzero finitely generated integral ideal of $S$ is invertible \cite[Definition 2.3]{GhalandarzadehNasehpourRazavi2017}. 
	
	\begin{theorem}
		
		\label{Dedekind1}
		
		The following conditions are equivalent for a semidomain $S$:
		
		\begin{enumerate}
			\item $S$ is a Noetherian Pr\"{u}fer semidomain.
			
			\item $\Frac(S)$ is an abelian group.
			
			\item $S$ is a Dedekind semidomain.
			
			\item $S$ is Noetherian and every 2-generated ideal of $S$ is invertible.
			
			\item Each nonzero integral ideal of $S$ is invertible.
		\end{enumerate}
		
		\begin{proof}
			
			$(1) \rightarrow (2)$: Since $S$ is a semidomain, $\Frac(S)$ is an abelian monoid. So, we only need to prove that each element of $\Frac(S)$ is invertible. Now, let $A$ be an arbitrary nonzero fractional ideal of $S$ and $d$ its common denominator. It is clear that $dA$ is an integral ideal of $S$. We observe that since $S$ is Noetherian, $dA$ is finitely generated and since $S$ is Pr\"{u}fer, $dA$ is invertible. Therefore, $A$ which is the multiplication of two invertible fractional ideals $dA$ and $(1/d)$, is itself invertible.
			
			The proof of the implications $(2) \rightarrow (3)$ and $(3) \rightarrow (4)$ is straightforward.
			
			$(4) \rightarrow (5)$: By Proposition \ref{Dedekindidentity}, the proof of this implication is straightforward.
			
			$(5) \rightarrow (1)$: By Proposition \ref{invertibleisfg}, $S$ is a Noetherian Pr\"{u}fer semidomain and the proof is complete.
		\end{proof}
	\end{theorem}
	
	\begin{example}
		Let $D$ be a Dedekind domain. By Theorem 3.7 in \cite{GhalandarzadehNasehpourRazavi2017}, the semiring of ideals $\Id(D)$ of $D$ is a Pr\"{u}fer semidomain. If $\Id(D)$ is also Noetherian, then $\Id(D)$ is a Dedekind semidomain. Note that the semiring $\Id(D)$ is proper, i.e., it is not a ring. For a more specific example, we assert that $(\Id(\mathbb Z), +,\cdot)$ is a principal ideal semidomain \cite[Remark 2.2]{GhalandarzadehNasehpourRazavi2017}. So, $\Id(\mathbb Z)$ is obviously a Dedekind semidomain. Note that the semiring $(\Id(\mathbb Z), +,\cdot)$ is isomorphic to the semiring $(\mathbb N_0, \gcd, \cdot)$. Finally, note that in Example \ref{DulinMosherEx1}, we show that the semidomain $\mathbb N_0$ under the usual addition and multiplication of non-negative numbers is not a Dedekind semidomain (in our sense defined in Definition \ref{Dedekindsemidomaindef}).
	\end{example}

	Using Theorem 2.9 in \cite{GhalandarzadehNasehpourRazavi2017} and this point that each Noetherian Pr\"{u}fer semidomain is Dedekind, we have the following corollary which gives some other equivalent definitions for Dedekind semidomains:
	
	\begin{theorem}
		
		\label{Dedekind2}
		
		Let $S$ be a Noetherian semidomain. Then the following statements are equivalent:
		
		\begin{enumerate}
			
			\item The semiring $S$ is a Dedekind semidomain,
			\item $\mathfrak{a}(\mathfrak{b}\cap \mathfrak{c})=\mathfrak{a}\mathfrak{b}\cap \mathfrak{a}\mathfrak{c}$ for all ideals $\mathfrak{a}$, $\mathfrak{b}$, and $\mathfrak{c}$ of $S$,
			\item $(\mathfrak{a}+\mathfrak{b})(\mathfrak{a}\cap \mathfrak{b})=\mathfrak{a}\mathfrak{b}$ for all ideals $\mathfrak{a}$ and $\mathfrak{b}$ of $S$,
			\item $[(\mathfrak{a}+\mathfrak{b}):\mathfrak{c}]=[\mathfrak{a}:\mathfrak{c}]+[\mathfrak{b}:\mathfrak{c}]$ for all ideals $\mathfrak{a}$, $\mathfrak{b}$, and $\mathfrak{c}$ of $S$,
			\item $[\mathfrak{a}:\mathfrak{b}]+[\mathfrak{b}:\mathfrak{a}]=S$ for all ideals $\mathfrak{a}$ and $\mathfrak{b}$ of $S$,
			\item $[\mathfrak{c}:\mathfrak{a}\cap \mathfrak{b}]=[\mathfrak{c}:\mathfrak{a}]+[\mathfrak{c}:\mathfrak{b}]$ for all ideals $\mathfrak{a}$, $\mathfrak{b}$, and $\mathfrak{c}$ of $S$.
			
		\end{enumerate}
	\end{theorem}
	
	I.S. Cohen in \cite{Cohen1950} proves that the domain $D$ is Dedekind if and only if it is Noetherian, and the lattice of its ideals is distributive, i.e., $\mathfrak{a} \cap (\mathfrak{b} + \mathfrak{c}) = (\mathfrak{a} \cap \mathfrak{b}) + (\mathfrak{a} \cap \mathfrak{c})$, for all ideals $\mathfrak{a}, \mathfrak{b},$ and $\mathfrak{c}$ in $D$. The following result is a special case of Lemma 2.4 in \cite{GhalandarzadehNasehpourRazavi2017}:
	
	\begin{proposition}
		
		\label{Dedekindisdistributive}
		
		If $S$ is a Dedekind semidomain, then the lattice of its ideals is distributive.
	\end{proposition}
	
	Theorem 8 in Cohen's paper \cite{Cohen1950} states that a Noetherian domain $D$ is Dedekind if and only if $D_{\mathfrak{m}}$ is discrete valuation ring, for each maximal ideal $\mathfrak{m}$ of $D$. We give a semiring version of this theorem in the following:
	
	\begin{theorem}
		
		\label{LocalizationofDedekind}
		
		For a semidomain $S$, the following statements are equivalent:
		
		\begin{enumerate}
			
			\item $S$ is a Dedekind semidomain.
			
			\item $S$ is Noetherian and for each maximal ideal $\mathfrak{m}$, $S_\mathfrak{m}$ is a discrete valuation semiring.
			
		\end{enumerate}
		
		\begin{proof}
			
			$(1) \rightarrow (2)$: By Theorem \ref{Dedekind1}, $S$ is Noetherian. Now let $\mathfrak{m}$ be a maximal ideal of $S$. It is clear that any ideal of $S_\mathfrak{m}$ is of the form $I_\mathfrak{m}$, where $I$ is an ideal of $S$. Now if $I$ is nonzero, then $I$ is invertible by hypothesis and so by Proposition 1.5 in \cite{GhalandarzadehNasehpourRazavi2017}, $I_\mathfrak{m}$ is principal. Also, note that $\mathfrak{m} S_\mathfrak{m} \neq (0)$. Now by Theorem 3.6 in \cite{Nasehpour2018(b)}, $S_\mathfrak{m}$ is a discrete valuation semiring.
			
			$(2) \rightarrow (1):$ Let $S$ be a Noetherian semiring and $I$ be a nonzero proper ideal of $S$. It is clear that $I$ is finitely generated. Let $\mathfrak{m}$ be an arbitrary maximal ideal of $S$. Since $S_\mathfrak{m}$ is a discrete valuation semiring, $I_\mathfrak{m}$ is nonzero and principal. Now by Theorem 1.8 in \cite{GhalandarzadehNasehpourRazavi2017}, $I$ is invertible and by using Theorem \ref{Dedekind1}, the proof is complete.
		\end{proof}
		
	\end{theorem}
	
	\begin{corollary}
		
		\label{locallyDedekind}
		
		If $S$ is a Dedekind semidomain, then so is $S_\mathfrak{m}$ for each maximal ideal $\mathfrak{m}$ of $S$. Moreover, if $S$ is a Noetherian semidomain, the converse is also true.
		
		\begin{proof}
			
			Let, for the moment, $S$ be a discrete valuation semiring. Then by Theorem 3.6 in \cite{Nasehpour2018(b)}, there exists a nonzero and nonunit element $t\in S$ such that any nonzero ideal $I$ of $S$ is of the form $I=(t^n)$ for $n\geq 0$. This clearly shows that any discrete valuation semiring (semidomain) is a Dedekind semidomain.
			
			Moreover, if $S$ is a Noetherian semidomain such that $S_\mathfrak{m}$ is Dedekind for each maximal ideal $\mathfrak{m}$ of $S$, then $S_\mathfrak{m}$ is a discrete valuation semiring. Now according to Theorem \ref{LocalizationofDedekind}, $S$ is Dedekind and the proof is complete.
		\end{proof}
		
	\end{corollary}
	
	\begin{remark}
		
		A semidomain $S$ is a Pr\"{u}fer semidomain if and only if $S_{\mathfrak{m}}$ is a valuation semidomain for each maximal ideal $\mathfrak{m}$ of $S$ \cite[Theorem 2.11]{GhalandarzadehNasehpourRazavi2017}. The Noetherian Pr\"{u}fer semidomains are exactly the Dedekind semidomains (Theorem \ref{Dedekind1}). It follows that if $S$ is a Dedekind semidomain then $S_{\mathfrak{m}}$ is a Noetherian valuation
		semidomain for each maximal ideal $\mathfrak{m}$ of $S$. Note that a local Noetherian valuation semidomain is a discrete valuation semiring (Theorem \ref{locallyDedekind}). By considering this discussion, we give the following definition:
		
		\begin{definition}
			We define a semidomain $S$ to be almost Dedekind if $S_\mathfrak{m}$ is Dedekind for each maximal ideal $\mathfrak{m}$ of $S$.
		\end{definition} 
		
		Clearly, Dedekind semidomains are almost Dedekind. Also, almost Dedekind semidomains are Pr\"{u}fer semidomains \cite[Theorem 2.11]{GhalandarzadehNasehpourRazavi2017}. Note that there are some commutative rings that they are almost Dedekind, i.e., locally Dedekind, but not Dedekind. Apparently, such examples first appeared in Nakano's paper \cite[p. 426]{Nakano1953}. More sophisticated examples can be found in Loper's note \cite{Loper2006}. So, the following question is quite natural:
		
	\end{remark}
	
	\begin{question}
		Is there any proper semiring such that it is locally a Dedekind semidomain while it is not itself a Dedekind semidomain?
	\end{question}
	
	\begin{definition}
		Given integral ideals $\mathfrak{a}$ and $\mathfrak{b}$ of any semidomain $S$, we say that $\mathfrak{a}$ divides $\mathfrak{b}$, denoted by $\mathfrak{a} | \mathfrak{b}$, if $\mathfrak{b} = \mathfrak{ac}$ for some integral ideal $\mathfrak{c}$.
	\end{definition}
	
	In any semidomain, if $\mathfrak{a}$ divides $\mathfrak{b}$ (i.e., $\mathfrak{b} = \mathfrak{ac}$ for some integral ideal $\mathfrak{c}$) then $\mathfrak{a}$ contains $\mathfrak{b}$. In each Dedekind semidomain, the converse is also true, as we show it in the following lemma:   
	
	\begin{lemma}
		
		\label{LemmaDivide}
		
		Let $\mathfrak{a}$ and $\mathfrak{b}$ be integral ideals in a Dedekind semidomain $S$. Then, the following statements are equivalent:
		
		\begin{enumerate}
			\item $\mathfrak{a}$ divides $\mathfrak{b}$,
			\item $\mathfrak{a}$ contains $\mathfrak{b}$,
			\item $\mathfrak{a}^{-1}$ is a subset of $\mathfrak{b}^{-1}$.
		\end{enumerate}
		
		\begin{proof}
			$(1) \rightarrow (2)$: Let $\mathfrak{a}$ divides $\mathfrak{b}$. So, there is an integral ideal $\mathfrak{c}$ of $S$ such that $\mathfrak{b} = \mathfrak{ac}$. Obviously,  $\mathfrak{ac} \subseteq \mathfrak{a} \cap \mathfrak{c} \subseteq \mathfrak{a}$. 
			
			$(2) \rightarrow (1)$: Let $\mathfrak{b} \subseteq \mathfrak{a}$. Since $\mathfrak{b}$ is invertible, there is a fractional ideal $\mathfrak{b}'$ such that $\mathfrak{b}\mathfrak{b}' = S$. Set $\mathfrak{c}= \mathfrak{a}\mathfrak{b}'$. It is clear that $\mathfrak{c}\subseteq S$. Now, we see that $\mathfrak{bc} = \mathfrak{b}\mathfrak{b}' \mathfrak{a}=\mathfrak{a}$.
			
			The proof of the equivalence of the statements (2) and (3) is straightforward.
		\end{proof}
	\end{lemma}
	
	In 1925, Krull \cite{Krull1925} investigated commutative rings in which the condition $\mathfrak{a} \subseteq \mathfrak{b}$ for ideals implies the existence of an ideal $\mathfrak{c}$ satisfying $\mathfrak{a} = \mathfrak{b}\mathfrak{c}$. He called them ``regular multiplication rings'' (in German, regul\"{a}res Multiplikationsring). Today these rings are called rather ``multiplication rings'' \cite[p. 71]{Gilmer1992}. Similarly, we give the following definition:
	
	\begin{definition}
		
		\label{multiplicationsemiring}
		
		We define a semiring $S$ to be multiplication if the condition $\mathfrak{a} \subseteq \mathfrak{b}$ for ideals $\mathfrak{a}$ and $\mathfrak{b}$ of $S$ implies the existence of an ideal $\mathfrak{c}$ satisfying $\mathfrak{a} = \mathfrak{bc}$.
	\end{definition}
	
	\begin{theorem}
		
		\label{Dedekindismultiplication}
		
		Let $S$ be a Noetherian semidomain. Then $S$ is multiplication if and only if $S$ is Dedekind.
		
		\begin{proof}
			Let $S$ be a Noetherian multiplication semidomain and let $\mathfrak{a}$ be an integral ideal of $S$. If $\mathfrak{a}$ is principal, then $\mathfrak{a}$ is invertible. If not, then there is a nonzero $c\in S$ such that $(c) \subset \mathfrak{a}$. So, there is an ideal $\mathfrak{b}$ of $S$ such that $(c) = \mathfrak{ab}$. Now, since $(c)$ is invertible, then $\mathfrak{a}$ is invertible. Therefore, by Theorem \ref{Dedekind1}, $S$ is a Dedekind semidomain. Conversely, let $S$ be a Dedekind semidomain. Clearly, by Lemma \ref{LemmaDivide}, $S$ is multiplication and this completes the proof.
		\end{proof}
		
	\end{theorem}
	
	\begin{definition}
		
		\label{weakmultiplicationsemiring}
		
		We define a semiring $S$ to be a weak multiplication semiring if $\mathfrak{a} \subseteq \mathfrak{p}$, where $\mathfrak{a}$ is an ideal of $S$ and $\mathfrak{p}$ is a prime ideal of $S$, implies that there exists an ideal $\mathfrak{c}$ of $S$ such that $\mathfrak{a} = \mathfrak{p}\mathfrak{c}$. 
	\end{definition}
	
	By Theorem \ref{multiplicationsemiring}, it is clear that each Dedekind semidomain is a weak multiplication semiring.
	
	\begin{theorem}
		
		\label{weakmultiplicationstrictlybetween}
		
		Let $S$ be a weak multiplication semiring. Then for its maximal ideals $\mathfrak{m}$, there are no ideals $\mathfrak{a}$ strictly between $\mathfrak{m}^2$ and $\mathfrak{m}$.
		
		\begin{proof} Let $\mathfrak{m}$ be a maximal ideal of $S$ and there is an ideal $\mathfrak{a}$ in $S$ such that $\mathfrak{m}^2 \subset \mathfrak{a} \subset \mathfrak{m}$. Since $S$ is weak multiplication, there is an ideal $\mathfrak{b}$ such that $\mathfrak{a}= \mathfrak{m} \mathfrak{b}$ and $\mathfrak{b} \nsubseteq \mathfrak{m}$. Let $x\in \mathfrak{b}-\mathfrak{m}$. Clearly, the principal ideal $(x)$ is a subset of $\mathfrak{b}$. So, $\mathfrak{m} \cdot (x)  \subseteq \mathfrak{m}\mathfrak{b} = \mathfrak{a}$. Finally, $\mathfrak{m}\cdot \big(\mathfrak{m} + (x)\big) = \mathfrak{m}^2 + \mathfrak{m} (x) \subseteq \mathfrak{a}$. Now, since $\mathfrak{m}$ is maximal, $\mathfrak{m}+(x) = S$, which implies $\mathfrak{m} \subseteq \mathfrak{a}$, a contradiction. This completes the proof. \end{proof}
		
	\end{theorem}
	
	It is a fact in ring theory that $D$ is a Dedekind domain if and only if $D$ is Noetherian, and for its maximal ideals $\mathfrak{m}$, there are no ideals $\mathfrak{a}$ strictly between $\mathfrak{m}^2$ and $\mathfrak{m}$ (see Theorem 8 in Cohen's paper \cite{Cohen1950}). We have the following result:
	
	\begin{corollary}
		
		\label{Dedekindstrictlybetween}
		
		Let $S$ be a Dedekind semidomain. Then for its maximal ideals $\mathfrak{m}$, there are no ideals $\mathfrak{a}$ strictly between $\mathfrak{m}^2$ and $\mathfrak{m}$.
	\end{corollary}
	
	\begin{question}
		Let $S$ be a Noetherian semidomain such that for its maximal ideals $\mathfrak{m}$, there are no ideals $\mathfrak{a}$ strictly between $\mathfrak{m}^2$ and $\mathfrak{m}$. Is $S$ a Dedekind semidomain?
	\end{question}
	
	Surprisingly, in ring theory, every weak multiplication ring is a multiplication ring \cite[Theorem 9.21]{LarsenMcCarthy1971}. Therefore, a Noetherian domain is Dedekind if and only if it is weak multiplication. For the author, it is not clear if the semiring version of this statement is true:
	
	\begin{question}
		Let $S$ be a Noetherian weak multiplication semidomain. Is $S$ Dedekind?
	\end{question}
	
	Let us recall that the Krull dimension $\dim S$ of a semiring $S$ is the supremum of the lengths of all chains of prime ideals in $S$ \cite[p. 578]{AlarconAnderson1994}.
	
	\begin{theorem}
		
		\label{Dedekindprimeismaximal}
		
		Each nonzero prime ideal of a Dedekind semidomain is maximal. In other words, the Krull dimension of a Dedekind semidomain is at most 1.

		\begin{proof}
			Let $S$ be a Dedekind semidomain and $\mathfrak{p}$ be a nonzero prime ideal of $S$. Since each proper ideal of $S$ is a subset of a maximal ideal of $S$, there is a maximal ideal $\mathfrak{m}$ of $S$ such that $\mathfrak{p} \subseteq \mathfrak{m}$. Therefore, by Lemma \ref{LemmaDivide}, there is an integral ideal $\mathfrak{c}$ of $S$ such that $\mathfrak{p}=\mathfrak{mc}$. Our claim is that $\mathfrak{p}=\mathfrak{m}$. On the contrary, suppose that $\mathfrak{p}\neq \mathfrak{m}$ and choose some $x\in \mathfrak{m} - \mathfrak{p}$. Vividly, $xc\in \mathfrak{p}$, for all $c\in \mathfrak{c}$ and since $\mathfrak{p}$ is prime, $\mathfrak{c}\subseteq \mathfrak{p}$. Using Lemma \ref{LemmaDivide}, this implies that $\mathfrak{c}=\mathfrak{p}\mathfrak{d}$, for some integral ideal $\mathfrak{d}$ of $S$. Now, we observe that \[\mathfrak{m}^{-1} = \mathfrak{p}^{-1}\mathfrak{c} = \mathfrak{d} \subseteq S. \]	Hence, $S \subseteq (\mathfrak{m}^{-1})^{-1} = \mathfrak{m}$, a contradiction and this proves the theorem.\end{proof}
	\end{theorem}

\begin{remark}
	
	\label{Krulldimatmost1}
	
	For Theorem \ref{Dedekindprimeismaximal}, we give the following alternative proof:
	
	\begin{proof}
	
	Let $\mathfrak{p}$ be a nonzero prime ideal of $S$. By Proposition 6.59 in \cite{Golan1999(b)}, $\mathfrak{p} \subseteq \mathfrak{m}$ for some maximal ideal $\mathfrak{m}$ of $S$. Since by Theorem \ref{LocalizationofDedekind}, $S_\mathfrak{m}$ is a discrete valuation semiring, by Theorem 3.6 in \cite{Nasehpour2018(b)}, the only nonzero prime ideal of $S_\mathfrak{m}$ is $\mathfrak{m} S_\mathfrak{m}$. Also, since $\mathfrak{p} S_\mathfrak{m} \neq (0)$, we have that $\mathfrak{p} S_\mathfrak{m} = \mathfrak{m} S_\mathfrak{m}$. Now we prove that $\mathfrak{p} = \mathfrak{m}$. Take $s\in \mathfrak{m}$. So $s/1 \in \mathfrak{m} S_\mathfrak{m} = \mathfrak{p} S_\mathfrak{m}$. This implies that $s/1 = t/u$, where $t\in \mathfrak{p}$ and $u \in S-\mathfrak{m}$. Now we have $us \in \mathfrak{p}$, while $u \notin \mathfrak{p}$. So $s\in \mathfrak{p}$ and the proof is complete. \end{proof}
\end{remark}
	
	\begin{example}
		
		\label{DulinMosherEx1}
		
		Let $(\mathbb N_0, +, \cdot)$ be the semiring of non-negative integer numbers under the usual addition and multiplication. Each subtractive ideal of $\mathbb N_0$ is of the form $m\mathbb N_0$, where $m$ is an arbitrary non-negative integer \cite[Proposition 6]{Noronha-Galvao1978}. Now, by the fundamental theorem of arithmetic, it is clear that each subtractive ideal of $\mathbb N_0$ is a product of subtractive prime ideals in $\mathbb N_0$. This means that the semiring $\mathbb N_0$ is Dedekind in the sense of Dulin and Mosher (check p. 87 in \cite{DulinMosher1972}). On the other hand, in $\mathbb N_0$,  $2\mathbb N_0$ is a prime ideal properly contained in the only maximal ideal $\mathbb N_0 - \{1\}$. So, by Theorem \ref{Dedekindprimeismaximal}, the semidomain $\mathbb N_0$ is not Dedekind in the sense of the current paper (defined in Definition \ref{Dedekindsemidomaindef}).
	\end{example}
	
	\begin{theorem}[The Unique Factorization Theorem for Ideals] 
		
		\label{UFTforIdeals}
		
		Let $A$ be a nonzero fractional ideal of a Dedekind semidomain $S$. Then, there are unique integers $v(\mathfrak{p},A)$, for $\mathfrak{p}\in \Spec(S)-\{0\}$, and almost all of which are 0, so that \[ A = \prod_{\mathfrak{p}\in \Spec(S)-\{0\}} \mathfrak{p}^{v(\mathfrak{p},A)} \label{UFTformula} \tag{UFT}.\]
		
		\begin{proof}
			
			Let $\mathfrak{a}$ be a proper nonzero integral ideal of $S$ such that it has no factorization of the form \eqref{UFTformula} mentioned above. Obviously, $\mathfrak{a}$ is not a maximal ideal of $S$. Thus there is some maximal ideal $\mathfrak{m}_1$ of $S$ which contains $\mathfrak{a}$ strictly. Using Lemma  \ref{LemmaDivide}, $\mathfrak{a}_1 = \mathfrak{a} \mathfrak{m}^{-1}_1$ is an integral ideal of $S$ such that $\mathfrak{a}\subset \mathfrak{a}_1$ (since $\mathfrak{m}_1 \subset S$). In succession, $\mathfrak{a}_1$ cannot be a maximal ideal of the semidomain $S$, for otherwise, we have found a factorization of the form \eqref{UFTformula} for $\mathfrak{a}$. So, in a similar way, we can find a maximal ideal $\mathfrak{m}_2$ which contains $\mathfrak{a}_1$ with $\mathfrak{a}_1 \subset \mathfrak{a}_1 \mathfrak{m}^{-1}_2 = \mathfrak{a}_2$. Continuing this way, we produce an infinite ascending chain of integral ideals in $S$, which is impossible since $S$ is Noetherian.
			
			Now, let $A$ be an arbitrary nonzero fractional ideal of $S$ and $d\in S$ a nonzero suitable common denominator for $A$ such that $dA \subseteq S$. So, $A = \mathfrak{b}\mathfrak{c}^{-1}$ for some nonzero integral ideals $\mathfrak{b}$ and $\mathfrak{c}$ of $S$. Finally, using factorization of the form \eqref{UFTformula} for the integral ideals $\mathfrak{b}$ and $\mathfrak{c}$, we obtain a factorization of the form \eqref{UFTformula} for $A$.
			
			Suppose that there are two factorizations of the form \eqref{UFTformula} for some nonzero fractional ideal $A$. By cross-multiplication (and cancellation of ideals if possible), we obtain an equality of the form \[\mathfrak{p}_1 \mathfrak{p}_2 \cdots \mathfrak{p}_k = \mathfrak{q}_1 \mathfrak{q}_2 \cdots \mathfrak{q}_l,\]where $\mathfrak{p}_1,\mathfrak{p}_2,\dots,\mathfrak{p}_k$ and $\mathfrak{q}_1,\mathfrak{q}_2,\dots,\mathfrak{q}_l$ are all nonzero prime ideals of $S$. Now, by Proposition \ref{FactorizationofInvertible}, the proof is complete.
		\end{proof}
		
	\end{theorem}
	
	We record three immediate corollaries for Theorem \ref{UFTforIdeals}.
	
	\begin{corollary}
		
		\label{UFTforintegralIdeals}
		
		Let $A$ be a nonzero fractional ideal of a Dedekind semidomain $S$. Then $A$ is integral if and only if $v(\mathfrak{p},A)\geq 0$, for all nonzero prime ideals $\mathfrak{p}$ of $S$.
	\end{corollary}
	
	\begin{corollary}
		Suppose that $\mathfrak{a}$ is an integral ideal of a Dedekind semidomain $S$.
		
		\begin{enumerate}
			\item There is a finite set of nonzero prime ideals $\mathfrak{p}$ which contain $\mathfrak{a}$, namely,
			those with $v(\mathfrak{a},\mathfrak{p}) \geq 1$.
			
			\item Let $\mathfrak{p}_1, \mathfrak{p}_2, \dots, \mathfrak{p}_k$ be the nonzero prime ideals containing $\mathfrak{a}$. Then the integral ideals which contain $\mathfrak{a}$ are those that can be written in the form \[B= \mathfrak{p}^{r_1}_1 \mathfrak{p}^{r_2}_2 \dots \mathfrak{p}^{r_k}_k \quad{with~ 0 \leq r_i \leq v(\mathfrak{p}_i, \mathfrak{a}) ~ for ~ all ~ i}.\] In particular, there is a finite set of integral ideals containing $\mathfrak{a}$.
		\end{enumerate}
		
	\end{corollary}

	\begin{corollary}
		Let $S$ be a Dedekind semidomain. Then the ideal group $\mathcal{I}_{\mathcal{F}}(S)$ of $S$ is a free abelian group with the collection of nonzero prime ideals of $S$ as free generators.
	\end{corollary}

Inspired by the proof of Theorem 12.6.8 in \cite{Steinberger1993}, we give the following result:

\begin{theorem}
	
	\label{ConverseUFTforIdeals}
	
	Let $S$ be a semidomain. Then the following statements are equivalent:
	
	\begin{enumerate}
		\item $S$ is a Dedekind semidomain.
		
		\item $S$ is a semidomain in which every nonzero integral ideal of $S$ is uniquely represented as a product of maximal ideals in $S$ (up to the order of the factors).
	\end{enumerate}

	\begin{proof}
		
		$(1) \rightarrow (2)$: By Theorem \ref{UFTforIdeals} and Corollary \ref{UFTforintegralIdeals}, every nonzero integral ideal of $S$ is uniquely represented as a product of prime ideals in $S$. On the other hand, by Theorem \ref{Dedekindprimeismaximal}, each nonzero prime is maximal.
		
		$(2) \rightarrow (1)$: Let $\mathfrak{m}$ be a nonzero maximal ideal of $S$ and choose $0\neq s\in \mathfrak{m}$. Let $ (s) = \mathfrak{m}_1 \cdots \mathfrak{m}_k$ be the factorization of the nonzero principal ideal $(s)$. So, $\mathfrak{m}_i \subseteq \mathfrak{m}$, for some $1\leq i \leq k$. Since $\mathfrak{m}_i$ and $\mathfrak{m}$ are both maximal ideals of $S$, $\mathfrak{m}_i = \mathfrak{m}$. Also, since $(s)$ is invertible, each factor of $(s)$ is also invertible which means that $\mathfrak{m}_i = \mathfrak{m}$ is invertible. Definitely, this implies that each nonzero ideal of $S$ is invertible and by Theorem \ref{Dedekind1}, $S$ is Dedekind and the proof is complete.   
	\end{proof}
	
\end{theorem}
	
	Let $S$ be a semidomain and $\mathcal{S}$ a nonempty subset of the set of all ideals $\Id(S)$ of $S$. We say the factorization law holds for $\mathcal{S}$ if each ideal in $\mathcal{S}$ may be factored into prime ideals in $\mathcal{S}$. The question arises if the factorization law for $\mathcal{S}$ causes the factorization law for $\Id(S)$. In this direction, one may ask the following natural question:
	
	Imagine the factorization law holds for the set of all subtractive ideals of $S$, i.e., $S$ is a Dedekind semidomain in the sense of Dulin and Mosher \cite{DulinMosher1972}. Will this cause the factorization law to hold for $\Id(S)$? The answer is negative, as the following example shows this: 
	
	\begin{example}
		
		\label{DulinMosherEx2}
		
		Consider the semiring $(\mathbb N_0, + , \cdot)$. In Example \ref{DulinMosherEx1}, we have already seen that the factorization law holds for the set of subtractive ideals of $\mathbb N_0$. On the other hand, each prime ideal of $\mathbb N_0$ is either of the form $p\mathbb N_0$, where $p$ is a prime number or is equal to the only maximal ideal of $\mathbb N_0$, i.e., $\mathbb N_0 - \{1\}$ \cite{AllenDale1975}. Our claim is that the ideal $I=\mathbb N_0 - \{1,2\}$ cannot be factored into primes of $\mathbb N_0$. On the contrary, let $I = P_1 P_2 \cdots P_m$, where $P_i$ is a prime in $\mathbb N_0$ for each $i$. The ideal $I=\mathbb N_0 - \{1,2\}$ is not prime, since $4\in I$, while $2\notin I$. So, the length of the factorization is at least 2, i.e., $m\geq 2$. Note that if $P$ is a prime ideal of $\mathbb N_0$, then for each element $x\in P$, we have $x\geq 2$. So, for each element $y \in P_1P_2 \cdots P_m$  with $m\geq 2$, we obtain that $y\geq 4$; while $3\in I$, a contradiction. This completes the proof of our claim. 
	\end{example}
	
	Mori worked on those commutative rings in which every principal ideal is a product of prime ideals \cite{Mori1938}. Such rings are now called $\pi$-rings \cite{Anderson1978}. In this direction, we give the following definition:
	
	\begin{definition}
		
		\label{pisemiring}
		
		Let $S$ be a semiring. We say $S$ is a $\pi$-semiring if each nonzero proper principal ideal of $S$ is a product of prime ideals of $S$.
		
	\end{definition}
	
	\begin{example}
		
		Consider the LaGrassa's semiring $L=\{0,u,1\}$, where $1+u = u+1 = u$, $1+1=1$, and $u+u = u\cdot u =u$ \cite{LaGrassa1995}. The only ideals of $L$ are $(0)$, $L$, and $(u)=\{0,u\}$. The ideal $\{0,u\}$ is maximal (and prime). Therefore, each nonzero proper principal ideal of $S$ is a product of prime ideals. Therefore, LaGrassa's semiring $L$ is an example for $\pi$-semirings. Note that $L$ is not a semidomain. 
		
	\end{example}

    Stephen McAdam, a student of Irving Kaplansky, proves that if $I$ is an ideal of a domain $D$ being maximal among all non-invertible ideals, then $I$ is prime \cite[Exercise 36, p. 44]{Kaplansky1974}. One of the corollaries of McAdam's statement is that if all nonzero primes in a nontrivial ring $R$ are invertible, then $R$ is a Dedekind domain \cite[Corollary 3.15]{LamReyes2008}.
	
	In the following, we find a general characterization for subtractive Dedekind semidomains:
	
	\begin{theorem}
		
		\label{SubtractiveDedekind1}
		
		Let $S$ be a subtractive semidomain. Then $S$ is a Dedekind semidomain if and only if the following two conditions hold:
		
		\begin{enumerate}
			
			\item \label{invertibleprimes} Every nonzero prime ideal in $S$ is invertible.
			
			\item \label{pisemidomain} $S$ is a $\pi$-semiring.
		\end{enumerate}
		
		\begin{proof}
			
			By Theorem \ref{Dedekind1} and Theorem \ref{UFTforIdeals}, the forward implication holds.
			
			Conversely, let the statements (\ref{invertibleprimes}) and (\ref{pisemidomain}) hold. At first, we prove that every nonzero prime ideal of $S$ is maximal, for if not, then there exists a prime ideal $\mathfrak{p}$ such that $(0) \subset \mathfrak{p} \subset \mathfrak{m} \subset S$, where $\mathfrak{m}$ is a suitable maximal ideal in $S$. Clearly, $\mathfrak{p}\mathfrak{m}^{-1}$ is an ideal of $S$ and $\mathfrak{p}=(\mathfrak{p}\mathfrak{m}^{-1})\mathfrak{m} \subseteq \mathfrak{p}\mathfrak{m}^{-1} \subseteq \mathfrak{p}$ and this implies that $\mathfrak{p}\mathfrak{m} =\mathfrak{p}$. Since $\mathfrak{p}$ is invertible, we obtain that $\mathfrak{m}=S$, a contradiction.  
			
			Since each nonzero prime ideal of $S$ is invertible, by Proposition \ref{invertibleisfg}, each prime ideal of $S$ is finitely generated and since $S$ is subtractive, by Proposition 7.17 in \cite{Golan1999(b)}, $S$ is Noetherian.
			
			Now, let $\mathfrak{a}$ be an arbitrary nonzero ideal of $S$. Since $S$ is Noetherian, we may write $\mathfrak{a} = (s_1, s_2, \dots, s_n)$ for some $s_i\in S$ such that at least one of them is nonzero. On the other hand, by assumption, each nonzero principal ideal of $S$ is a product of maximal ideals (nonzero primes) in $S$. Therefore, by Proposition \ref{IdealOperationsCoprime}, \[\mathfrak{a} = (s_1) + (s_2) + \dots + (s_n)\] is also a product of maximal ideal in $S$. Now, we see that each nonzero ideal of $S$ is invertible and by Theorem \ref{Dedekind1} $S$ is Dedekind and this is what we were supposed to prove.\end{proof}
		
	\end{theorem}
	
	\begin{remark}
		
		Cohen shows that if every nonzero ideal in the domain $D$ is a product of prime ideals, then $D$ is a Dedekind domain (see Theorem 6 in \cite{Cohen1950}). Then with the help of this statement, he proves that if every nonzero prime ideal in the domain $D$ is invertible, then $D$ is a Dedekind domain. Also, note that in some resources an integral domain in which every nonzero ideal is uniquely represented as the product of a finite number of prime ideals is called a Dedekind domain (see for example p. 294 in Matsumura's book \cite{Matsumura1980}). Based on this, we give the following questions:
		
		\begin{questions}
			
			\begin{enumerate}
				\item Let $S$ be a semidomain in which every nonzero ideal is a product of prime ideals in $S$. Is $S$ a Dedekind semidomain?
				
				\item Let $S$ be a semidomain in which every nonzero ideal is uniquely represented as the product of a finite number of prime ideals in $S$. Is $S$ a Dedekind semidomain?
			\end{enumerate}
			
		\end{questions}
	\end{remark}
	
	Inspired by E. Matlis \cite[p. 258]{Matlis1970}, we define $\mu_{*}(S)$ to be the supremum of the minimum number of generators required to generate an ideal of $S$. It is a famous fact in ring theory that if $D$ is a Dedekind domain, then $\mu_{*}(D) \leq 2$. We prove its semiring version in the following:
	
	\begin{theorem}
		
		\label{SubtractiveDedekind2}
		
		If $S$ is a subtractive Dedekind semidomain, then $\mu_{*}(S) \leq 2$.
		
		\begin{proof}
			
			Let $I$ be a nonzero ideal of $S$. Therefore, there is a nonzero element $a\in I$ such that $(a) \subseteq I$. So, by Theorem \ref{UFTforIdeals}, we know that \[(a) = \mathfrak{m}_1^{k_1} \cdots \mathfrak{m}_n^{k_n} \quad \text{and} \quad I= \mathfrak{m}_1^{l_1} \cdots \mathfrak{m}_n^{l_n}, \] where $\mathfrak{m}_1, \dots, \mathfrak{m}_n$ are distinct maximal ideals of $S$ and $l_i \leq k_i$ are natural numbers for all $1\leq i \leq n$. By uniqueness property in Theorem \ref{UFTforIdeals}, for each $1\leq i \leq n$, there is a $b_i \in S$ such that \[b_i \in \mathfrak{m}_1^{l_1 + 1} \cdots \mathfrak{m}_i^{l_i} \cdots \mathfrak{m}_n^{l_n + 1} - \mathfrak{m}_1^{l_1 + 1} \cdots \mathfrak{m}_i^{l_i + 1} \cdots \mathfrak{m}_n^{l_n + 1} \subseteq I.\] On the other hand, $b_i \in \mathfrak{m}_j^{l_j +1}$ for each $j\neq i$, while $b_i \notin \mathfrak{m}_i^{l_i +1}$, for if not, \[b_i \in \mathfrak{m}_i^{l_i +1} \cap (\mathfrak{m}_1^{l_1 + 1} \cdots \mathfrak{m}_i^{l_i} \cdots \mathfrak{m}_n^{l_n + 1}) = \mathfrak{m}_1^{l_1 + 1} \cdots \mathfrak{m}_i^{l_i + 1} \cdots \mathfrak{m}_n^{l_n + 1},\] which is in contradiction with the choice of $b_i$.
			
			Define $b:= b_1 + \dots + b_n$. Since $\mathfrak{m}_i^{l_i +1}$ is subtractive, $b \notin \mathfrak{m}_i^{l_i +1}$ for each $1\leq i \leq n$, while definitely $b\in I$
			
			Since $(a,b) = (a)+(b)$, the prime factorization of $(a, b)$ contains at most the maximal ideals occurring
			in $(a)$, so we obtain the following \[ (a,b) = \mathfrak{m}_1^{t_1} \cdots \mathfrak{m}_n^{t_n}.\]
			
			Note that
			
			\begin{enumerate}[($\alpha$)]
				\item $l_i \leq t_i$, since $(a,b) \subseteq I$;
				
				\item $t_i \leq l_i$, since $(a,b) \nsubseteq \mathfrak{m}_i^{l_i +1}$ ($b \notin \mathfrak{m}_i^{l_i +1}$).
			\end{enumerate}
			So, $m_i=l_i$, which means that $I = (a,b)$. Hence, $\mu_{*}(S) \leq 2$ and the proof is complete.
		\end{proof}
		
	\end{theorem}
	
	\begin{theorem}
		
		\label{FiniteSpecPISD}
		
		A subtractive Dedekind semidomain with only finitely many prime ideals is a PISD.
		
		\begin{proof}
			
			Let $\mathfrak{a}$ be an arbitrary nonzero ideal of a subtractive Dedekind semidomain $S$ with the only distinct nonzero prime ideals $\{\mathfrak{p}_i\}_i^n$. Define $\mathfrak{b}:= \mathfrak{p}_1 \cdots \mathfrak{p}_n$. 
			
			Our claim is that there is a nonzero ideal $\mathfrak{a}^{\prime}$ such that $\mathfrak{a} \mathfrak{a}^{\prime}$ is a principal ideal and $\mathfrak{a}^{\prime}$ and $\mathfrak{b}$ are coprime. We prove our claim as follows:
			
			Define $\mathfrak{q}_i := \mathfrak{p}_1 \cdots \mathfrak{p}_{i-1} \mathfrak{p}_{i+1} \cdots \mathfrak{p}_n$, for each $i$. By Theorem \ref{UFTforIdeals}, $\mathfrak{a} \mathfrak{q}_i \supset \mathfrak{a} \mathfrak{b}$. For each $i$, choose $a_i \in \mathfrak{a} \mathfrak{q}_i - \mathfrak{a} \mathfrak{b}$ and let $a:= a_1 + \dots + a_n$. Clearly, $a\in \mathfrak{a}$, while similar to the proof of Theorem \ref{SubtractiveDedekind2}, $a_i \notin \mathfrak{a} \mathfrak{p}_i$. Also, for each $j\neq i$, we have $a_j \in \mathfrak{a} \mathfrak{q}_j \subseteq \mathfrak{a} \mathfrak{p}_i$. 
			
			Since each ideal of $S$ is subtractive, $a\notin \mathfrak{a} \mathfrak{p}_i$, for each $i$.
			
			On the other hand, $(a)\in \mathfrak{a}$. Since $S$ is Dedekind, by Theorem \ref{Dedekindismultiplication}, there is an ideal $\mathfrak{a}^{\prime}$ such that $\mathfrak{a} \mathfrak{a}^{\prime} = (a)$. In order to complete the proof of our claim, we show that $\mathfrak{p}_i$ cannot divide $\mathfrak{a}^{\prime}$, for each $i$. Because if this is so for some $i$, then $\mathfrak{a}^{\prime} = \mathfrak{p}_i \mathfrak{a}_0$ for some nonzero ideal $\mathfrak{a}_0$. This implies that $(a) = \mathfrak{a} \mathfrak{p}_i \mathfrak{a}_0$ which leads us to the contradiction $a\in \mathfrak{a} \mathfrak{p}_i$.
			
			So, the ideals $\mathfrak{a}^{\prime}$ and $\mathfrak{b}$ are coprime and this implies $\mathfrak{a}^{\prime}=S$. Now, $\mathfrak{a} = \mathfrak{a} S = \mathfrak{a} \mathfrak{a}^{\prime} = (a)$. This means that each ideal of $S$ is principal and this completes the proof.
		\end{proof}
		
	\end{theorem}

\begin{remark}(Examples of subtractive Dedekind semidomains) \label{ExamplessubtractiveDedekindsemidomains} Let $D$ be a Dedekind semidomain. Then $\Id(D)$ is a subtractive Pr\"{u}fer semidomain \cite[Theorem 3.7]{GhalandarzadehNasehpourRazavi2017}. Now, if $\Id(D)$ is Noetherian, then $\Id(D)$ is a subtractive Dedekind semidomain.
\end{remark}
	
	\section{$M$-Cancellation Ideals and $M$-Gaussian Semirings}\label{sec:Mcancel}
	
	Let us recall a nonzero ideal $\mathfrak{a}$ of a semiring $S$ is called a cancellation ideal if $\mathfrak{ab}=\mathfrak{ac}$ implies $\mathfrak{b}=\mathfrak{c}$ for all ideals $\mathfrak{b}$ and $\mathfrak{c}$ of $S$ \cite{LaGrassa1995}. Also, if $R$ is a commutative ring with a nonzero identity and $M$ is a unital $R$-module, then an ideal $\mathfrak{a}$ of $R$ is defined to be $M$-cancellation if for all $R$-submodules $P$ and $Q$ of $M$, $\mathfrak{a}P = \mathfrak{a}Q$ implies $P=Q$ \cite[Definition 2.1]{NasehpourYassemi2000}. Similarly, we give the following definition:
	
	\begin{definition}
		
		\label{Mcancellationdef}
		
		Let $S$ be a semiring and $M$ an $S$-semimodule. We define an ideal $\mathfrak{a}$ of $S$ to be $M$-cancellation if for all $S$-subsemimodules $P$ and $Q$ of $M$, $\mathfrak{a}P = \mathfrak{a}Q$ implies $P=Q$.
	\end{definition}
	
	\begin{proposition}
		
		\label{invertibleiscancellation}
		
		Let $M$ be an $S$-semimodule. Every integral invertible ideal of a semiring $S$ is $M$-cancellation.
		
		\begin{proof} Straightforward. \end{proof}
	\end{proposition}
	
	Let us recall that an $R$-module $M$  is an Auslander module if $r\in R$ is not a zero-divisor on $M$, then $r$ is not a zero-divisor on $R$, or equivalently, if $Z_R(R) \subseteq Z_R(M)$ \cite{NasehpourAuslander}. On the other hand, if $R$ is a ring, $M$ an $R$-module, and $Q$ the total ring of fractions of $R$, then $M$ is torsion-free if the natural map $M \rightarrow M \otimes Q$ is injective \cite[p. 19]{BrunsHerzog1998}. It is straightforward to see that $M$ is a torsion-free $R$-module if and only if $Z_R(M) \subseteq Z_R(R).$
	
	\begin{proposition}
		
		\label{canMcan}
		
		Let $R$ be a commutative ring with a nonzero identity and $M$ a unital $R$-module. Then the following statements hold:
		
		\begin{enumerate}
			\item If $M$ is not an Auslander module (i.e. there is an element $r\in R$ such that $r$ is a zero-divisor on $R$ while it is not a zero-divisor on $M$), then there is a principal ideal in $R$ such that it is an $M$-cancellation, while not a cancellation ideal of $R$.
			
			\item If $M$ is not a torsion-free module (i.e. there is an element $r\in R$ such that $r$ is a zero-divisor on $M$ while it is not a zero-divisor on $R$), then there is a principal ideal in $R$ such that it is a cancellation, while not an $M$-cancellation ideal of $R$.
			
			\item The set of principal cancellation ideals is the same as the set of principal $M$-cancellation ideals in $R$ if and only if the set of zero-divisor elements on $R$ is equal to the set of zero-divisor elements on $M$, i.e. $Z_R(R) = Z_R(M)$.
		\end{enumerate}
		
		\begin{proof} Straightforward. \end{proof}
		
	\end{proposition}
	
	\begin{remark}
		To obtain some examples satisfying the conditions in Proposition \ref{canMcan}, see Remark 3.3 and Proposition 3.4 in \cite{NasehpourAuslander}. 
	\end{remark}
	
	Let us recall that if $S$ is a semiring and $X$ is an indeterminate, the set of all polynomials over the semiring $S$, denoted by $S[X]$, is the set of all formal forms $a_0 + a_1 X + \cdots + a_n X^n$, where $a_0, a_1, \ldots, a_n \in S$. Similar to ring theory, $S[X]$ is a semiring under the usual addition and multiplication of polynomials. In the same way, if $M$ is an $S$-semimodule, one can consider $M[X]$ as an $S[X]$-semimodule under the standard addition and scalar product.
	
	For any $g \in M[X]$, the content of $g$, denoted by $c(g)$, is defined to be the $S$-subsemimodule of $M$ generated by the coefficients of $g$. \cite{Nasehpour2016}.
	
	By considering the semimodule version of the Dedekind-Mertens Lemma (see \cite{Northcott1950} and Theorem 2 in \cite{Nasehpour2016}) and the concept of Gaussian semirings (check Definition 7 in \cite{Nasehpour2016}), we define the following concept:
	
	\begin{definition}
		
		\label{MGaussiandef}
		
		Let $M$ be an $S$-semimodule. We say a semiring $S$ is $M$-Gaussian if $c(fg) = c(f)c(g)$ for all $f\in S[X]$ and $g\in M[X]$. 
	\end{definition}
	
	\begin{theorem}
		
		\label{MGaussian}
		
		Let $S$ be a Pr\"{u}fer (in particular, a Dedekind) semidomain and $M$ be a subtractive $S$-semimodule. Then $S$ is $M$-Gaussian.
		
		\begin{proof}
			
			Let $f\in S[X]$ and $g\in M[X]$. By Theorem 2 in \cite{Nasehpour2016}, there is a non-negative integer $n$ such that $c(f)^{n+1} c(g)= c(f)^n c(fg)$. Since $S$ is Pr\"{u}fer, each finitely generated ideal of $S$ is invertible and so, $M$-cancellation. So, by canceling $c(f)^n$, we get the equality $c(fg) = c(f)c(g)$ and the proof is complete. \end{proof}
		
	\end{theorem}
	
	\begin{corollary}
		Let $S$ be a subtractive Pr\"{u}fer semidomain. Then $S$ is Gaussian.
	\end{corollary}
	
	Let us recall that if $B$ is an $S$-semialgebra, the content of an element $f\in B$, denoted by $c(f)$, is defined to be the following ideal: \[c(f) = \bigcap \{ I \colon I \text{~is an ideal of~} S \text{~and~} f \in IB \}.\] By Definition 30 in \cite{Nasehpour2016}, $B$ is a content $S$-semialgebra if $S$ is a subsemiring of $B$ and the following conditions hold:
	
	\begin{enumerate}
		
		\item $f\in c(f)B$ for all $f\in B$;
		
		\item $c(sf) = sc(f)$ for all $s\in S$ and $f\in B$ and $c(1)= S$;
		
		\item (Dedekind-Mertens content formula) For all $f,g \in B$ there exists an $m\in \mathbb N_0$ such that $c(f)^{m+1} c(g) = c(f)^m c(fg)$.
		
	\end{enumerate}
	
	Let $R$ be a ring. An $R$-algebra $B$ is called Gaussian if $c(fg)=c(f)c(g)$ for all elements $f,g \in B$ (see Definition 1 in \cite{NasehpourGaussian}). Similarly, we define Gaussian semialgebra:
	
	\begin{definition}
		We define an $S$-semialgebra $B$ to be Gaussian if it is a content semialgebra and $c(fg)=c(f)c(g)$ for all elements $f,g \in B$.
	\end{definition}
	
	\begin{proposition}
		Let $S$ be a Pr\"ufer (in particular, a Dedekind) semidomain. If $B$ is a content $S$-semialgebra, then it is a Gaussian $S$-semialgebra.
		
		\begin{proof} Let $f,g\in B$ be arbitrary. If one of these elements is zero, say $f$, then $c(f)=(0)$ and so, there is nothing to prove. Now, take both $f$ and $g$ to be nonzero. By definition, there is a non-negative $n$ such that $c(f)^{n+1} c(g) = c(f)^n c(fg)$. Since $B$ is a content $S$-semialgebra, by Proposition 23 in \cite{Nasehpour2016}, $c(f) \neq (0)$ is finitely generated and so invertible (since $S$ is Pr\"ufer). By canceling $c(f)^n$, we get $c(fg)=c(f)c(g)$ and this completes the proof. \end{proof}
		
	\end{proposition}
	
	\section*{Acknowledgments} The author is supported in part by the Department of Engineering Science at the Golpayegan University of Technology and wishes to thank the department for supplying all necessary facilities in pursuing this research. The author is also grateful to Professor Dara Moazzami for his help and encouragements.

\end{document}